\documentclass[graybox]{svmult}
\usepackage{amsmath}
\usepackage{amsfonts}
\usepackage[english]{babel}
\usepackage[latin1]{inputenc}
\usepackage[T1]{fontenc}
\usepackage{ae}
\usepackage{url}
\allowdisplaybreaks[4]
\usepackage{color}
\usepackage{tikz}
\usetikzlibrary{matrix}
\usepackage{cite}           
\usepackage{mathptmx}
\usepackage{helvet}
\usepackage{courier}
\usepackage{makeidx}
\usepackage{graphicx} 
\usepackage{multicol} 
\usepackage[bottom]{footmisc} 
\newcommand{\et}[2]{\eta_{#1}^{#2}}
\makeindex
\begin{document}
\title*{Eta quotients\index{Eta quotients} and Rademacher sums\index{Rademacher sums}}
\author{Kevin Acres and David Broadhurst}
\institute{Kevin Acres, Monash University, Victoria 3800, Australia,\\
\email{research@research-systems.com},\\
David Broadhurst, Open University, Milton Keynes MK76AA, UK,\\
\email{David.Broadhurst@open.ac.uk}}
\maketitle
\abstract{Eta quotients on $\Gamma_0(6)$ yield evaluations
of sunrise integrals at 2, 3, 4 and 6 loops.
At 2 and 3 loops, they provide modular parametrizations\index{modular parametrizations}
of inhomogeneous differential equations whose solutions
are readily obtained by expanding in the nome $q$. 
Atkin-Lehner transformations\index{Atkin-Lehner transformations}
that permute cusps ensure fast convergence for all
external momenta.  At 4 and 6 loops, on-shell integrals
are periods of modular forms of weights 4 and 6 given by 
Eichler integrals\index{Eichler integrals} of eta quotients.
Weakly holomorphic eta quotients determine quasi-periods\index{quasi-periods}.
A Rademacher sum formula is given for Fourier
coefficients of an eta quotient that is a Hauptmodul\index{Hauptmodul} for 
$\Gamma_0(6)$ and its generalization is found for
all levels with genus 0, namely for
$N = 1,2, 3, 4, 5, 6, 7, 8, 9, 10, 12, 13, 16, 18, 25$.
There are elliptic obstructions 
at $N = 11, 14, 15, 17, 19, 20, 21, 24, 27, 32, 36, 49,$
with genus 1. We surmount these,
finding explicit formulas for Fourier coefficients of eta quotients in thousands of cases.
We show how to handle the levels $N=22, 23, 26, 28, 29, 31, 37, 50$, with genus 2,
and the levels $N=30,33,34,35,39,40,41,43,45,48,64$, with genus 3.
We also solve examples with genera $4,5,6,7,8,13$.} 

\section{Introduction}

Elliptic obstructions to the evaluation of massive Feynman diagrams
were recognized and surmounted more than 50 years ago by Afaf Sabry~\cite{S}.
They occur in two-loop two-point integrals when three massive
particles appear in an intermediate state~\cite{Bold}. The simplest
example is the two-loop sunrise diagram with unit masses
in two space-time dimensions,
whose study was revolutionized in 2013, when Spencer Bloch and
Pierre Vanhove~\cite{BV} showed how to parametrize and solve its second
order differential equation using eta quotients on $\Gamma_0(6)$.

Their solution was particularly bold, since they expand in a nome $q$
that is small near the physical threshold where the external energy
$w$ is close to 3. Thus they achieve fast convergence near the
branchpoint that frustrates other methods. The price to pay
is that convergence is slow near any of the other three
cusps of $\Gamma_0(6)$, which occur at $w=0,1,\infty$.
We shall show how to use Atkin-Lehner transformations
of eta quotients to expand about those cusps, 
achieving optimal efficiency.

Bloch, Kerr and Vanhove~\cite{BKV} conquered the corresponding
three-loop problem, also using eta quotients on $\Gamma_0(6)$,
thanks to the remarkable circumstance,
noted more than 40 years ago by Geoffrey Joyce~\cite{Joyce},
that a transformation of variables relates solutions of the 
relevant homogeneous third-order differential equation
to products of solutions of the second-order equation at two loops.
Joyce's observation was made in the context of the physics
of condensed matter. The relevance of his work on the
diamond lattice to Feynman integrals
was decoded in~\cite{BBBG}. We shall use
an Atkin-Lehner transformation to achieve
optimal efficiency at three loops.

The role of $\Gamma_0(6)$ does not end at three loops.
It is of the essence for the on-shell problems at 4 and 6 loops,
where the relevant Bessel moments turn out to be
Eichler integrals of eta quotients that are cusp forms
of level 6 with modular weights 4 and 6, respectively.
We shall review key results, which were until
recently only conjectures~\cite{mzvmf,filskm,BS,BM},
tested to many thousands of digits.
For an account of how they were proved~\cite{Z1,Z2,Z3,Z4},
see the lucid review by Yajun Zhou~\cite{Z5}.

It is notable that this connection between number theory 
and Feynman integrals persists in the
real world of four-dimensional space-time.
The four-loop radiative corrections to the
magnetic moment of the electron in 
quantum electrodynamics,
evaluated with breath-taking skill by Stefano Laporta~\cite{Laporta},
contain a pair of Bessel moments~\cite{Z4}
that are Eichler integrals. We conclude Section~2 
with results that indicate that one of these
is a quasi-period, in the sense of Francis Brown~\cite{Brown}.
Moreover we conjecturally identify
quasi-periods at 6 loops.
 
Section 3 concerns a searching question raised by Johannes  Bl\"umlein
at a recent conference held at the Hausdorff Centre for Mathematics, in Bonn.
Is there a closed formula for the Fourier coefficients of
the Hauptmodul of $\Gamma_0(6)$, of the type that
Petersson~\cite{P} and Rademacher~\cite{R,R1,K} found
for Klein's $j$-invariant?
We conjecturally answer in the affirmative, by giving
a formula that serves this purpose for all levels with genus 0.
Moreover we are able to extend its use to higher genera.

\section{Eta quotients in quantum field theory}

Broadhurst, Fleischer and Tarasov~\cite{BFT} gave the differential  
equation for the two-loop unit-mass sunrise integral in an arbitrary number $D$
of space-time dimensions. At $D=2$, this integral is a Bessel moment~\cite{BBBG}
\begin{equation}
I(w^2)=4\int_0^\infty I_0(wx)K_0^3(x)x{\rm d}x,
\end{equation}
where $w$ is the external energy, which enters the Bessel function $I_0(wx)$ via Fourier transformation.
The Bessel function $K_0(x)$ is cubed, since three particles of unit mass connect the two vertices.
Bloch and Vanhove~\cite{BV} found a very neat modular parametrization of the differential equation at $D=2$,
which we here write as  
\begin{gather} 
-\left(q\frac{{\rm d}}{{\rm d}q}\right)^2\frac{I(w^2)}{6f}=g,\quad w=\frac{3\et22\et34}{\et14\et62},
\quad f=\frac{\et16\et6{}}{\et23\et32},\\
g=\frac{\et25\et34\et6{}}{\et14}=\frac{\et39}{\et13}+\frac{\et69}{\et23}
=\sum_{n>0}\frac{n^2(q^n-q^{5n})}{1-q^{6n}},
\quad\eta_n=q^{n/24}\prod_{k>0}(1-q^{nk}),
\end{gather}
with eta quotients determining the energy $w$, the integrating factor $f$, which is an elliptic integral
determining the discontinuity across the cut for $w>3$, and the inhomogeneous term
$g$. Two easy integrations of the Lambert series for $g$ then yield
\begin{equation}
\frac{I(w^2)}{f}=\frac{\pi\log(-1/q)}{\sqrt3}
-3\sum_{n>0}\frac{\chi_6(n)}{n^2}\frac{1+q^n}{1-q^n}\label{ex1}
\end{equation}
with $\chi_6(n)=\pm1$ for $n=\pm1$ mod 6 and $\chi_6(n)=0$, otherwise. This solution
is determined by the discontinuity across the cut and the
finiteness of $I(1)=\pi^2/4$~\cite{BBBG}.
In summary: after dividing $I(w^2)$ by the modular form $f$, with weight 1 and level 6, %
we obtain solution~(\ref{ex1}) by two integrations of the weight 3 modular form $g$ %
with respect to $z$, where $q=\exp(2\pi{\rm i}z)$. Such integrals of modular forms are %
referred to as Eichler integrals. %

We remark that modular parametrizations of differential equations were used in~\cite{FB}, %
to elucidate proofs of rationality of zeta values, and in~\cite{KZ}, for problems in statistical physics.   %

\subsection{Atkin-Lehner transformations of eta quotients}

Now set $q=\exp(2\pi{\rm i}z)$ with $\Im z> 0$ and consider the transformations~\cite{CZ}
\begin{equation}
z\longmapsto z_2=\frac{2z-1}{6z-2},\quad 
z\longmapsto z_3=\frac{3z-2}{6z-3},\quad
z\longmapsto z_6=\frac{-1}{6z},\label{AL}
\end{equation}
which permute the cusps at $z=0,\frac12,\frac13,\infty$. Then, with $q_k=\exp(2\pi{\rm i}z_k)$,
\begin{gather}
-\left(q_k\frac{{\rm d}}{{\rm d}q_k}\right)^2\frac{I(w^2)}{6f_k(z_k)}=g_k(z_k),\label{alt}\\
f_2(z)=\frac{\et26\et3{}}{\et13\et62},\quad
f_3(z)=\frac{\et2{}\et36}{\et12\et63},\quad
f_6(z)=\frac{\et1{}\et66}{\et22\et33},\\
g_2(z)=\frac{\et15\et3{}\et64}{\et24},\quad
g_3(z)=\frac{\et14\et2{}\et65}{\et34},\quad
g_6(z)=\frac{\et1{}\et24\et35}{\et64}.
\end{gather}

From the alternative differential equations~(\ref{alt}), we obtain alternative expansions:
\begin{eqnarray}
\frac{I(w^2)}{f_2(z_2)}&=&\sum_{n=1}^\infty\frac{3\chi_3(n)}{n^2}\frac{(1-q_2^n)^2}{1+q_2^{2n}}
=I(0)-\sum_{n=1}^\infty\frac{6\chi_3(n)}{n^2}\frac{q_2^n}{1+q_2^{2n}},\label{ex2}\\
\frac{I(w^2)}{f_3(z_3)}&=&\sum_{n=1}^\infty\frac{2\chi_2(n)}{n^2}\frac{(1-q_3^n)^3}{1-q_3^{3n}}
=I(1)-\sum_{n=1}^\infty\frac{6\chi_2(n)}{n^2}\frac{q_3^n}{1+q_3^n+q_3^{2n}},\label{ex3}\\
\frac{I(w^2)}{f_6(z_6)}&=&-3\log^2(-q_6)
+\sum_{n=1}^\infty\frac{6}{n^2}\frac{q_6^n}{1-q_6^n+q_6^{2n}}\label{ex6},
\end{eqnarray}
with $\chi_2(n)=0,1$, for $n=0,1$ mod 2, and $\chi_3(n)=-1,0,1$, for $n=-1,0,1$ mod 3.

Then for any real value of $w^2$ there is an optimal choice of nome in which
to expand, which may be determined as follows. Let
\begin{equation}
w_1^2=w^2,\quad w_2^2=\frac{w^2-9}{w^2-1},\quad w_3^2=\frac{9}{w^2},\quad w_6^2=\frac{9}{w_2^2}. 
\end{equation}
For $w^2\in[-3,\,9-6\sqrt2]$ set $k=2$, else for $w^2\in[9-6\sqrt2,\,3]$ set $k=3$, 
else for $w^2\in[3,\,9+6\sqrt2]$ set $k=1$, else set $k=6$. Then compute
$w_k\in[\sqrt3,\sqrt3+\sqrt6]$ and obtain the optimal nome $q_k=Q(w_k)$ from
\begin{equation}
Q(x)=\exp\left(\frac{-\pi\,{\rm agm}(1,\sqrt{r})}{{\rm agm}(1,\sqrt{1-r})}\right),\quad r=\frac{16x}{(x+3)(x-1)^3},
\end{equation}
by the lightning-fast process of the arithmetic-geometric mean. 
This results in a small real nome  $q_k\in[-\exp(-\pi/\sqrt3),\,\exp(-\pi\sqrt{2/3})]$ and hence 
$|q_k|<0.16304$. If $k=2$, use~(\ref{ex2}); if $k=3$, use~(\ref{ex3}); if $k=6$, use~(\ref{ex6});
if $k=1$ use $q=q_1$ in~(\ref{ex1}) and extract a Clausen value from
\begin{equation}
\sum_{n>0}\frac{\chi_6(n)}{n^2}\frac{1+q^n}{1-q^n}=
C_2+\sum_{n>0}\frac{\chi_6(n)}{n^2}\frac{2q^n}{1-q^n},\quad
C_2=\frac{5\,{\rm Cl}_2(\pi/3)}{\sqrt{27}}=\frac{5I(0)}{12}.
\end{equation}

The authors of~\cite{Mainz} expand in $q_{M}=-q_2$, thereby encountering $\eta_4$ and $\eta_{12}$ in
\begin{equation}
f_M(z)=f_2\left(z+\frac12\right)=\frac{\et13\et43\et6{}}{\et23\et3{}\et{12}{}},\quad
g_M(z)=g_2\left(z+\frac12\right)=-\frac{\et2{11}\et67}{\et15\et3{}\et45\et{12}{}}.
\end{equation}
Since they expand about the cusp at $w=0$, they inevitably face issues 
of slow convergence near the cusps at $w=1,3,\infty$. Moreover they had to address
delicate questions of analytic continuation at the on-shell point $w=1$. 
Our procedure of invariably expanding about the {\em nearest} cusp avoids 
all such problems.

Such use of Atkin-Lehner transformations to achieve efficient expansions in small nomes is well known %
to mathematicians who compute with modular forms~\cite{GP}. %
We recommend exploitation of the transformations~(\ref{AL}) to physicists who encounter problems %
that involve the congruence subgroup $\Gamma_0(6)$. For example, the authors of~\cite{Abl} encountered, at 3 loops,
a second-order equation with complicated coefficients and powers of $\log(x)$ in the inhomogeneous term.
Their homogeneous equation has a hypergeometric solution
\begin{equation}
H(x)=\frac{(x^2-1)^2}{9(x^2+3)}
\sum_{n=0}^\infty\frac{(4/3)_n(5/3)_n}{n!(n+1)!}\left(\frac{x^2(x^2-9)^2}{(x^2+3)^3}\right)^{n+1}
\end{equation}
where $(a)_n=\Gamma(a+n)/\Gamma(a)$ is the Pochhammer symbol. We obtained
\begin{equation}H\left(3\frac{\et12\et64}{\et24\et32}\right)=
\frac12\left(\frac{\et1{14}\et6{10}}{\et2{22}\et32}
+\frac{\et16\et64}{\et2{12}\et32}\left(\frac{\et14\et68}{\et28\et34}+\frac13\right)
q\frac{{\rm d}}{{\rm d}q}\right)\frac{\et26\et3{}}{\et13\et62}
\end{equation}
as a modular parametrization of the homogeneous solution, where the derivative
with respect to $q$ results from a complete integral of the second kind.
It would be interesting to investigate whether an Atkin-Lehner transformation
may be used to avoid a singularity that was encountered at $x=\frac13$
at intermediate stages of the work in~\cite{Abl}.
 
For our next advertisement of the virtue of Atkin-Lehner transformation, we turn
to the three-loop equal-mass sunrise integral.
Bailey, Borwein, Broadhurst and Glasser~\cite{BBBG} developed the expansion in $t$ of
\begin{equation}
J(t)=8\int_0^\infty I_0(\sqrt{t}x)K_0^4(x)x{\rm d}x=7\zeta(3)+O(t).
\end{equation}
A neat and novel $q$-expansion comes from exploiting a transformation~\cite{BBBG,Joyce}
from $w^2$, at two loops, to $t=10-w^2-9/w^2$, at three loops.
Then we obtain the modular parametrization
\begin{gather}
t=10-w^2-\frac{9}{w^2}=-64\left(\frac{\eta_2\eta_6}{\eta_1\eta_3}\right)^6,\\
\left(q\frac{{\rm d}}{{\rm d}q}\right)^3\frac{J(t)}{(wf/3)^2}
=24h,\quad h=\frac{\et2{16}}{\et18}-9\frac{\et6{16}}{\et38}
=\sum_{n>0}\frac{n^3(q^n-8q^{3n}+q^{5n})}{1-q^{6n}},\\
\frac{J(t)}{(wf/3)^2}=J(0)+24\sum_{n>0}\frac{\phi(n)}{n^3}\frac{q^n}{1-q^n},\label{Jex}
\end{gather}
with $\phi(n)=0,1,0,-8,0,1$, for $n=0,1,2,3,4,5$ mod 6.
The {\tt Pari-GP} procedure
\begin{verbatim}
z(t)={local(x=2/(sqrt(4-t)+sqrt(16-t)),
a=sqrt((1-x)^3*(1+3*x)));
I/2*agm(a,4*x*sqrt(x))/agm(a,sqrt((1+x)^3*(1-3*x)));}
\end{verbatim}
returns the correct value of $z$ for the nome $q=\exp(2\pi{\rm i}z)$, for all real $t$.
Expansion~(\ref{Jex}) is highly efficient for $t\in[-8,8]$, where $|q|\le\exp(-\sqrt2\pi/3)<0.22742.$
For the rest of the real $t$-axis, we exploit the involution $z\longmapsto z_6=-1/(6z)$, which gives
$t\longmapsto t_6=64/t$, with fixed points at $t=\pm8$. For $t_6\in[-8,8]$ we use
\begin{gather}
\left(q_6\frac{{\rm d}}{{\rm d}q_6}\right)^3\frac{J(t)}{(wf_6(z_6))^2}=-24h_6(z_6),\\
h_6(z)=-t_6h=1+2h-30\sum_{n>0}\frac{n^3(q^{2n}+q^{4n}-8q^{6n})}{1-q^{6n}},\\
\frac{J(t)}{(wf_6(z_6))^2}=-4\log^3(q_6)+24\sum_{n>0}\frac{15\phi(n+3)-\phi(n)}{n^3}\frac{1+q_6^n}{1-q_6^n},
\label{sol6}
\end{gather}
in agreement with the result proved by Bloch, Kerr and Vanhove~\cite{BKV}. Extracting
\begin{equation}
C_3=\sum_{n>0}\frac{15\phi(n+3)-\phi(n)}{n^3}=\frac{2\zeta(3)}{3}
\end{equation}
we achieve a highly efficient expansion in $q_6=\exp(-\pi{\rm i}/(3z))$ for $64/t\in[-8,8]$, with a strong check
of consistency with~(\ref{Jex}) in the neighbourhoods of $t=\pm8$, where both
expansions work well.

\subsection{Eichler integrals of eta quotients for on-shell sunrise integrals}

On-shell sunrise integrals lead us to consider Bessel moments of the form
\begin{equation}
M(a,b,c)=\int_0^\infty I_0^a(x)K_0^b(x)x^c{\rm d}x.
\end{equation}
For $L>3$, the off-shell $L$-loop integral 
\begin{equation}S_L(t)=\int_0^\infty\frac{{\rm d}x_1}{x_1}\ldots\int_0^\infty\frac{{\rm d}x_L}{x_L}
\frac{1}{(1+\sum_{j=1}^L x_j)(1+\sum_{k=1}^L 1/x_k)-t}\label{SL}
\end{equation}
has not yielded to the methods given above for $L=2,3$. By contrast, the on-shell values
$S_L(1)=2^LM(1,L+1,1)$ with $L+2$ Bessel functions yield Eichler integrals of 
cusp forms of weights 4 and 6 on $\Gamma_0(6)$
at $L=4$ and $L=6$ loops,
namely integrals of the form $\int_0^\infty f({\rm i}y)y^{s-1}{\rm d}y$, %
where $f(z)$ is a cusp form with weight $L$ and $s$ is a integer with $L>s>0$. %
First we consider the situation at $L=3$ loops, where a modular form of weight 3 occurs. 

At 3 loops, with 5 Bessel functions, the on-shell problem is solved by the weight 3 level 15 cusp form
\begin{equation}
f_{3,15}(z)=(\eta_3\eta_5)^3+(\eta_1\eta_{15})^3=\sum_{n>0}A_5(n)q^n=
-\frac{f_{3,15}(-1/(15z))}{(-15)^{3/2}z^3}
\end{equation}
with complex multiplication in ${\bf Q}(\sqrt{-15})$. 
If the Kronecker symbol $\left(\frac{p}{15}\right)=\left(\frac{p}{3}\right)\left(\frac{p}{5}\right)$ is negative, 
for prime $p$, then $A_5(p)=0.$ 
For $\Re s>2$, there is a convergent L-series
\begin{equation}
L_5(s)=\sum_{n>0}\frac{A_5(n)}{n^s}=\prod_{p}\frac{1}{1-A_5(p)p^{-s}+\left(\frac{p}{15}\right)p^{2-2s}}
\end{equation}
whose analytic continuation is provided by the Eichler integral
\begin{equation}
L_5(s)=\frac{(2\pi)^s}{\Gamma(s)}\int_0^\infty f_{3,15}({\rm i}y)y^{s-1}{\rm d}y
\end{equation}
with critical values
\begin{equation}
L_5(1)=\frac{5}{\pi^2}M(1,4,1),\quad
L_5(2)=\frac43M(2,3,1).
\end{equation}

At 4 loops, with 6 Bessel functions, the on-shell problem is solved by the weight 4 level 6 cusp form
\begin{equation}
f_{4,6}(z)=(\eta_1\eta_2\eta_3\eta_6)^2=\sum_{n>0}A_6(n)q^n
=\frac{f_{4,6}(-1/(6z))}{6^2z^4}.
\end{equation}
For $\Re s>5/2$, there is a convergent L-series
\begin{equation}
L_6(s)=\sum_{n>0}\frac{A_6(n)}{n^s}=
\frac{1}{1+2^{1-s}}\frac{1}{1+3^{1-s}}
\prod_{p>3}\frac{1}{1-A_6(p)p^{-s}+p^{3-2s}}
\end{equation}
whose analytic continuation is provided by the Eichler integral
\begin{equation}
L_6(s)=\frac{(2\pi)^s}{\Gamma(s)}\int_0^\infty f_{4,6}({\rm i}y)y^{s-1}{\rm d}y
\end{equation}
with critical values
\begin{gather}
L_6(2)=\frac{2}{\pi^2}M(1,5,1)=\frac{2}{3}M(3,3,1),\\
L_6(1)=\frac{2}{\pi^2}M(2,4,1)=\frac{3}{\pi^2}L_6(3).
\end{gather}

At 6 loops, with 8 Bessel functions, the on-shell problem is solved by the weight 6 level 6 cusp form
\begin{equation}
f_{6,6}(z)=\frac{\et29\et39}{\et13\et63}+\frac{\et19\et69}{\et23\et33}
=\sum_{n>0}A_8(n)q^n
=-\frac{f_{6,6}(-1/(6z))}{6^3z^6}.
\end{equation}
For $\Re s>7/2$, there is a convergent L-series
\begin{equation}
L_8(s)=\sum_{n>0}\frac{A_8(n)}{n^s}=
\frac{1}{1-2^{2-s}}\frac{1}{1+3^{2-s}}
\prod_{p>3}\frac{1}{1-A_8(p)p^{-s}+p^{5-2s}}
\end{equation}
whose analytic continuation is provided by the Eichler integral
\begin{equation}
L_8(s)=\frac{(2\pi)^s}{\Gamma(s)}\int_0^\infty f_{6,6}({\rm i}y)y^{s-1}{\rm d}y
\end{equation}
with critical values
\begin{gather}
L_8(4)=\frac{4}{9\pi^2}M(1,7,1)=\frac{4}{9}M(3,5,1)=\frac{\pi^2}{9}L_8(2),\\
L_8(5)=\frac{4}{27}M(2,6,1)=\frac{2\pi^2}{21}M(4,4,1)=\frac{2\pi^2}{21}L_8(3)=\frac{\pi^4}{54}L_8(1).
\end{gather}

\subsection{Eichler integrals for quasi-periods at level 6}

In~\cite{Brown} Francis Brown associated a pair of periods and a pair of
quasi-periods to the weight 12 level 1 modular form $\Delta(z)=\et1{24}$. The periods are
a pair of Eichler integrals that determine critical values of the L-series at odd
and even integers. No concrete integrals were given for the quasi-periods.
Rather it was asserted that numerical values may be obtained by an
undeclared regularization of integrals of a weakly holomorphic
modular form $\Delta^\prime(z)=1/q+O(q^2)$. 

In the case of the level 6 modular forms that yield 4-loop and 6-loop Feynman integrals
the situation is cleaner, since the periods
are Eichler integrals of eta quotients, $f_{4,6}$ and $f_{6,6}$, with 4 cusps.
Thus we may hope to find weakly holomorphic modular forms,
$g_{4,6}$ and $g_{6,6}$, 
that yield quasi-periods as well defined Eichler integrals
with a base-point at a cusp free of singularities.
A test is provided by the condition that a $2\times2$
determinant formed from a pair of periods and a pair of quasi-periods
should be an algebraic multiple of a power of $\pi$, as is trivially ensured
for modular forms of weight 2 by Legendre's relation
between pairs of complete elliptic integrals of first and second kind.

At 4 loops, we achieved this with Eichler integrals
\begin{gather}
\quad\frac{D_2}{2}=\frac{M(1,5,1)}{\pi^4}=\frac{4M(1,5,3)}{\pi^4} + \frac{5E_2}{18}\\
\frac{3D_1}{5}=\frac{M(2,4,1)}{\pi^3}=\frac{4M(2,4,3)}{\pi^3} + \frac{E_1}{3}\\
\left[\begin{array}{c}D_s\\E_s\end{array}\right]=-\int_{1/\sqrt3}^\infty
\left[\begin{array}{c}f_{4,6}\left(\frac{1+{\rm i}y}{2}\right)\\g_{4,6}\left(\frac{1+{\rm i}y}{2}\right)\end{array}\right]
y^{s-1}{\rm d}y,\\
g_{4,6}(z)=\frac{(w^2-3)^2(w^4+9)}{8w^4}f_{4,6}(z)=5q + 102q^2 + 945q^3 +O(q^4),\\
D_1E_2-D_2E_1=\frac{1}{24\pi^3}.
\end{gather}

At 6 loops, it is conjecturally achieved by
\begin{gather}
\qquad\qquad{\rm det}\left[\begin{array}{lr}
M(1,7,1)&32M(1,7,3)-64M(1,7,5)\\
M(2,6,1)&32M(2,6,3)-64M(2,6,5)\end{array}\right]=\frac{5\pi^6}{192},\\
\qquad\;\frac{F_2}{4}=\frac{M(1,7,1)}{\pi^6}\;\stackrel{?}{=}\;\frac{32M(1,7,3)-64M(1,7,5)}{\pi^6}+\frac{35G_2}{108},\\
\quad\frac{9F_1}{28}=\frac{M(2,6,1)}{\pi^5}\;\stackrel{?}{=}\;\frac{32M(2,6,3)-64M(2,6,5)}{\pi^5}+\frac{5G_1}{12},\\
\left[\begin{array}{c}F_s\\G_s\end{array}\right]=-\int_{1/\sqrt3}^\infty
\left[\begin{array}{c}f_{6,6}\left(\frac{1+{\rm i}y}{2}\right)\\g_{6,6}\left(\frac{1+{\rm i}y}{2}\right)\end{array}\right]
(3y^2-1)y^{s-1}{\rm d}y,\\
g_{6,6}(z)=\frac{(w^2-3)^4}{16w^4}f_{6,6}(z)=q + 36q^2 + 567q^3+5264q^4+O(q^5),\\
F_1G_2-F_2G_1\;\stackrel{?}{=}\;\frac{1}{4\pi^5},
\end{gather}
where the question marks indicate unproven discoveries, checked to thousands of digits of numerical precision.

\section{Rademacher sums for Fourier coefficients of eta quotients}

For positive integers $N$, $M$ and $n$, we define the Rademacher sums
\begin{equation}
R_{N,M}(n)=\sum_{c>0,~{\rm gcd}(c,N)=1}
\frac{2\pi I_{1}(4\pi\sqrt{nM/N}/c)}{\sqrt{nN/M}c}K(c,N,M,n)\label{RNM}
\end{equation}
as sums of Bessel functions multiplied by Kloosterman sums
\begin{equation}
K(c,N,M,n)=\sum_{r\in[1,c],~{\rm gcd}(r,c)=1}
\left.\exp\left(\frac{2\pi{\rm i}(Mr-ns)}{c}\right)\right|_{Nrs~=~1~{\rm mod}~c}\label{Kc}.
\end{equation}
In~(\ref{RNM}) the sum is over all positive integers $c$ coprime to $N$.
In~(\ref{Kc}) the sum is over the integers $r\in[1,c]$ coprime to $c$
and $s\in[1,c]$ is the inverse of $Nr$ modulo $c$. It follows from these
definitions that $R_{N,M}(n)/M=R_{N,n}(M)/n$ and that
$R_{N,M}(n)=R_{dN,dM}(n)$ for every positive integer $d$ that divides $N$.  

\subsection{Genus 0}

We found that $R_{N,1}(n)/R_{N,1}(1)$ is the coefficient of $q^n$ in an eta quotient
$T_N/B_N$ defining an OEIS sequence in the genus 0 cases of Table~\ref{table:1},
where the eta quotients agree with the canonical Hauptmoduln in~\cite[Table~8]{Maier}.

\begin{table}[h!]
\centering
\begin{tabular}{|c|c|c|c|c|} 
\hline
$N$&$R_{N,1}(1)$&$T_N$&$B_N$&OEIS\\
\hline
2&  4096&$\et2{24}$&$\et1{24}$&A014103\\
3&   729&$\et3{12}$&$\et1{12}$&A121590\\
4&   256&$\et48$&$\et18$&A092877\\
5&   125&$\et56$&$\et16$&A121591\\
6&    72&$\et2{}\et65$&$\et15\et3{}$&A128638\\
7&    49&$\et74$&$\et14$&A121593\\
8&    32&$\et22\et84$&$\et14\et42$&A107035\\
9&    27&$\et93$&$\et13$&A121589\\
10&  20&$\et2{}\et{10}3$&$\et13\et5{}$&A095846\\
12&  12&$\et22\et3{}\et{12}3$&$\et13\et4{}\et62$&A187100\\
13&  13&$\et{13}2$&$\et12$&A121597\\
16&   8&$ \et2{}\et{16}2$&$\et12\et8{}$&A123655\\
18&   6&$\et2{}\et3{}\et{18}2$&$\et12\et6{}\et9{}$&A128129\\
25&   5&$\et{25}{}$&$\et1{}$&A092885\\
\hline
\end{tabular}
\caption{Eta quotients $T_N/B_N$ of genus 0 with Fourier coefficients $R_{N,1}(n)/R_{N,1}(1)$}
\label{table:1}
\end{table}

For $n>0$, Rademacher~\cite{R} obtained $R_{1,1}(n)$ as the coefficient of $q^n$ of
\begin{equation}
j(z)=\frac{1}{\et1{24}}\left(1+240\sum_{n>0}\frac{n^3q^n}{1-q^n}\right)^3=\frac{1}{q}+744+196884q+O(q^2)\label{j}
\end{equation}
which is invariant under $z\mapsto(az+b)/(cz+d)$ with integers satisfying $ad-bc=1$.

We remark that analytic continuation of (\ref{RNM}) to $n=0$ gives $R_{1,1}(0)=24$, which differs from %
the constant term 744 in~(\ref{j}). Our work concerns only the values of $R_{N,M}(n)$ for integers $n>0$. %

The congruence subgroup $\Gamma_0(N)$ is the group of M\"obius transformations with $N|c$.
The Hauptmodul 
\begin{equation}\frac{\et2{}\et65}{\et15\et3{}}=
q + 5q^2 + 19q^3 + 61q^4 + 174q^5 + 455q^6 + 1112q^7+\ldots
\end{equation}
of $\Gamma_0(6)$
has a Fourier coefficient $R_{6,1}(n)/72={\rm A}128638(n)$, 
which we are now able to evaluate at large $n$. We found that this Fourier coefficient is odd %
if and only if the core (i.e. the square-free part) of $n$ is a divisor of 6. %
We determined the probably prime values of ${\rm A}128638(n)$ for $n\in[1,900000000]$ and %
found these occur at surprisingly few values of $n$, namely these: %
2, 3, 4, 9, 32, 48, 324, 578, 864, 121032, 940896, 11723776,
88360000, 180848704, 198443569.

We remark that
${\rm A}128638(900000000)$, with 66832 decimal digits, would be rather 
hard to compute in the absence of a Rademacher-sum formula.

\subsection{Further examples of integer sequences}

We found several integer sequences of the form $R_{N,M}(n)/D$,
with ${\rm gcd}(N,M)=1$, $N>M>1$ and integer $D$, as for example in Table~\ref{table:2}.

\begin{table}[h!]
\centering
\begin{tabular}{|c|c|c|l|} 
\hline
$N$&$M$&$D$&sequence\\
\hline
3 & 2& $3^7$& 8, 339, 6552, 82796, 790896, 6171606, 41232064, 243306300,\ldots\\
4 & 3& $2^{10}$ & 33, 1800, 42412, 633024, 7003278, 62405984, 471069624, 3114275328,\ldots\\
5 & 2& $5^3$& 12, 197, 1824, 12426, 68780, 327819, 1391472, 5383270, 19289244,\ldots\\
6 & 5& 432& 145, 10085, 286435, 5004925, 63619086, 642751655, 5445694040,\ldots\\
7 & 2& $7^2$& 8, 81, 504, 2476, 10248, 37590, 125328, 387384, 1123992, 3092369,\ldots\\
8 & 3& $2^7$& 9, 132, 1132, 7200, 37566, 169648, 685368, 2532096, 8688909,\ldots\\
9 & 2& $3^4$& 2, 15, 72, 287, 984, 3051, 8704, 23286, 58968, 142677, 331728,\ldots\\
10 & 3& 80 &6, 63, 418, 2139, 9216, 35004, 120594, 384147, 1146842, 3241083,\ldots\\
11& 8&$11^2$& 234, 11950, 266994, 3812019, 40551362, 348772038, 2548265460,\ldots\\
12& 5& 72& 25, 435, 4255, 30255, 174126, 859305, 3766760, 15014775, 55334545,\ldots\\
13& 2& 13& 4, 21, 72, 222, 600, 1509, 3536, 7902, 16860, 34740, 69264, 134412, \ldots\\
14& 5& 56& 17, 229, 1852, 11213, 55998, 243084, 946991, 3382221, 11242933,\ldots\\
15& 7& 45& 67, 1398, 15919, 128386, 826187, 4509396, 21688133, 94244610,\ldots\\
16& 3&$2^5$& 3, 18, 76, 264, 810, 2264, 5880, 14400, 33583, 75132, 162180, 339296,\ldots\\
17& 3& 17& 5, 26, 107, 352, 1045, 2814, 7091, 16842, 38225, 83260, 175329,\ldots\\
18& 5& 36& 10, 95, 580, 2770, 11226, 40340, 132080, 401255, 1145740, 3104412,\ldots\\
19& 2& 19& 1, 4, 10, 25, 55, 116, 229, 440, 809, 1455, 2541, 4354, 7300, 12050,\ldots\\
\hline
\end{tabular}
\caption{Examples of integer sequences $R_{N,M}(n)/D$}
\label{table:2}
\end{table}

We identified some of the generating functions, as follows
\begin{eqnarray}
\sum_{n>0}R_{3,2}(n)q^n&=&\frac{3^7\et3{12}}{\et1{12}}\left(8+\frac{3^5\et3{12}}{\et1{12}}\right)\\
\sum_{n>0}R_{5,2}(n)q^n&=&\frac{5^3\et56}{\et16}\left(12+\frac{5^3\et56}{\et16}\right)\\
\sum_{n>0}R_{7,2}(n)q^n&=&\frac{7^2\et74}{\et14}\left(8+\frac{7^2\et74}{\et14}\right)\\
\sum_{n>0}R_{9,2}(n)q^n&=&\frac{3^4\et93}{\et13}\left(2+\frac{3^2\et93}{\et13}\right)\\
\sum_{n>0}R_{13,2}(n)q^n&=&\frac{13\et{13}2}{\et12}\left(4+\frac{13\et{13}2}{\et12}\right)\\
\sum_{n>0}R_{16,3}(n)q^n&=&8\left(\frac{\et2{18}}{\et1{12}\et46}-1\right).
\end{eqnarray}

Moreover,
\begin{eqnarray}
\sum_{n>0}R_{18,2}(n)q^n&=&4\left(\frac{\et26\et32}{\et16\et62}-1\right)\\
\sum_{n>0}R_{27,2}(n)q^n&=&3\left(\frac{\et34}{\et13\et9{}}-1\right)\\
\sum_{n>0}R_{32,3}(n)q^n&=&4\left(\frac{\et22\et44}{\et14\et82}-1\right)\\
\sum_{n>0}R_{36,5}(n)q^n&=&6\left(\frac{\et25\et35}{\et17\et63}-1\right)\\
\sum_{n>0}R_{48,7}(n)q^n&=&6\left(\frac{\et28\et34}{\et18\et42\et62}-1\right)\\
\sum_{n>0}R_{64,3}(n)q^n&=&2\left(\frac{\et2{}\et42}{\et12\et8{}}-1\right)\\
\sum_{n>0}R_{64,7}(n)q^n&=&4\left(\frac{\et27}{\et16\et8{}}-1\right).
\end{eqnarray}

When $N>1$ has genus 0, $R_{N,M}(n)$ is an integer sequence generated
by a polynomial  of degree $M$ in the eta quotient that generates $R_{N,1}(n)$.
Thus, for example, we may compute the 
coefficient of $q^n$ in $(\eta_{25}/\eta_1)^P$ from a linear combination of Rademacher-type formulas 
for $R_{25,M}(n)$ with $M\in[1,P]$, using polynomials in 
$g=\sum_{n>0}R_{25,1}(n)q^n=5\eta_{25}/\eta_1$ as follows:
\begin{eqnarray}
\sum_{n>0}R_{25,2}(n)q^n&=&2g+g^2\\
\sum_{n>0}R_{25,3}(n)q^n&=&6g+3g^2+g^3\\
\sum_{n>0}R_{25,4}(n)q^n&=&12g+10g^2+4g^3+g^4\\
\sum_{n>0}R_{25,5}(n)q^n&=&25g+25g^2+15g^3+5g^4+g^5=\frac{5^3\et56}{\et16}\\
\sum_{n>0}R_{25,6}(n)q^n&=&42g+60g^2+44g^3+21g^4+6g^5+g^6\\
\sum_{n>0}R_{25,7}(n)q^n&=&77g+126g^2+119g^3+70g^4+28g^5+7g^6+g^7\\
\sum_{n>0}R_{25,8}(n)q^n&=&120g+260g^2+288g^3+210g^4+104g^5+36g^6+8g^7+g^8.\quad{}
\end{eqnarray}

\subsection{Genus 1}

The genus $g_0(N)$ of $\Gamma_0(N)$ is computed in {\tt Pari-GP} by a procedure
\begin{verbatim}
g0(N)={local(f=factor(N),t=vector(4,k,1),p,r,n);
for(k=1,matsize(f)[1],p=f[k,1];r=f[k,2];n=p^r;
t[1]*=n*(1+1/p);t[2]*=if(n==2,1,if(p%4==1,2));
t[3]*=if(n==3,1,if(p%3==1,2));
t[4]*=if(r%2,2*p^((r-1)/2),(p+1)*p^(r/2-1)));
1+t[1]/12-t[2]/4-t[3]/3-t[4]/2;}
\end{verbatim}
that combines 4 multiplicative functions~\cite{EZ,M,SZ}.
\begin{table}[h!]
\centering
\begin{tabular}{|l|r|} 
\hline
$N$&primes $M<1000$\\
\hline
11& 19, 29, 199, 569, 809\\
14& 5, 11, 23, 71, 101, 263, 503\\
15& 7, 23, 31, 79, 167, 431, 479, 983\\
17& 3, 11, 47, 359, 967\\
19& 2, 23, 257, 449, 509, 521, 641\\
20& 11, 131, 251, 491, 599\\
21& 23, 31, 47, 71, 127, 367, 383, 743\\
24& 7, 47, 191, 383, 439\\
27& $M=2$ mod 3\\
32& $M=3$ mod 4\\
36& $M=5$ mod 6\\
49& $M=3,5,6$ mod 7\\
\hline
\end{tabular}
\caption{Primes $M$ such that $R_{N,M}(n)$ is an integer sequence for $N$ with genus 1}
\label{table:3}
\end{table}

We conjecture that only when $N$ has genus 0 is $R_{N,1}(n)$ an integer sequence.
To deal with genus 1, we introduced the additional parameter $M$ into $R_{N,M}(n)$.
For each level $N$ with genus 1, we specify in Table~\ref{table:3} the 
prime values of $M<1000$, coprime to $N$, for which $R_{N,M}(n)$ 
is an integer sequence.

At genus 1, the criterion for whether $R_{N,M}(n)$ forms
an integer sequence is provided by the Fourier expansion of the unique 
weight 2 cusp form of level $N$, which we denote by
$f_N=\sum_{M>0}C_{N,M}q^M$. Specifically,
\begin{eqnarray}
f_{11}&=&(\eta_1\eta_{11})^2\\
f_{14}&=&\eta_1\eta_2\eta_7\eta_{14}\\
f_{15}&=&\eta_1\eta_3\eta_5\eta_{15}\\
f_{17}&=&q-q^2-q^4-2q^5+4q^7+3q^8-3q^9+2q^{10}-2q^{13}+\ldots\\
f_{19}&=&q-2q^3-2q^4+3q^5-q^7+q^9+3q^{11}+4q^{12}-4q^{13}\ldots\\
f_{20}&=&(\eta_2\eta_{10})^2\\
f_{21}&=&q-q^2+q^3-q^4-2q^5-q^6-q^7+3q^8+q^9+2q^{10}+4q^{11}-q^{12}+\ldots\quad{}\\
f_{24}&=&\eta_2\eta_4\eta_6\eta_{12}\\
f_{27}&=&(\eta_3\eta_9)^2\\
f_{32}&=&(\eta_4\eta_8)^2\\
f_{36}&=&\et64\\
f_{49}&=&q+q^2-q^4-3q^8-3q^9+4q^{11}-q^{16}-3q^{18}+4q^{22}\ldots
\end{eqnarray}
with explicit formulas for $N=17,19,21,49$ given below in~(\ref{lat1},\ref{lat2},\ref{lat3}).

For $N$ with genus 1, we found that $R_{N,M}(n)$ is an integer sequence 
if and only if $C_{N,M}=0$. Moreover 
\begin{equation}
\overline{R}_{N,M}(n)=R_{N,M}(n)-C_{N,M}R_{N,1}(n)
\end{equation}
is always an integer sequence, with $\overline{R}_{N,1}(n)=0$, by construction.  

With $G_{N,M}=\sum_{n>0}\overline{R}_{N,M}(n)q^n$, we found at $N=21$ that
\begin{eqnarray}
7\left(\frac{         \et3{} \et73}          {\et13\et{21}{}}     -1\right)&=&G_{21,2},\label{e21a}\\
3^3\left(\frac{   \et37  \et7{}}              {\et17\et{21}{} }-1\right)&=&G_{21,4}+G_{21,3}+2G_{21,2},\label{e21b}\\
7^2\left(\frac{    \et32 \et76}               {\et16\et{21}2  }-1\right) &=&G_{21,4}+2G_{21,3}+5G_{21,2},\label{e21c}\\
\frac{3^37    \et34 \et{21}2}                 {\et16               }&=&G_{21,4}-2G_{21,3}-G_{21,2},\\
\frac{3^37^2\et3{}\et{21}5}                 {\et15\et7{}       }&=&G_{21,4}-5G_{21,3}+5G_{21,2},\\
\frac{3^37    \et36 \et72}                     {\et18               }&=&G_{21,5}-2G_{21,2},\\
\frac{3^37^2\et33 \et7{}\et{21}3}         {\et17               }&=&G_{21,5}-3G_{21,4}+4G_{21,2}.
\end{eqnarray}

\begin{table}[h!]
\centering
\begin{tabular}{|c|r|r|c|} 
\hline
$N$ & $C_{N,2}$ & $C_{N,3}$ & $E_N$\\ %
\hline
11 & -2 & -1 & $Y^2 + 6XY + 121Y  = X^3 + 38X^2 + 363X$ \\  %
14 & -1  & -2 & $Y^2 + 3XY +  56Y  = X^3 + 25X^2 + 168X$ \\    %
15 & -1 &  -1 & $ Y^2 + 3XY +  45Y  = X^3 + 23X^2 + 135X$ \\ %
17 & -1 &   0 & $Y^2 + 3XY +  34Y  = X^3 + 18X^2 +  85X$ \\   %
19 & 0  &  -2 & $Y^2         +  19Y  =  X^3 + 16X^2 +  76X$ \\   %
20 & 0  &  -2 & $Y^2         +  20Y  = X^3 + 13X^2 +  60X$ \\      %
21 & -1 &  1 & $Y^2 + 3XY +  21Y  = X^3 + 13X^2 +  42X$ \\   %
24 &  0 & -1 & $Y^2         +  12Y  = X^3 + 11X^2 +  36X$ \\     %
27 &  0 &  0 & $Y^2         +   9Y  = X^3 +  9X^2 +  27X$ \\      %
32 & 0  &  0 & $Y^2         +   8Y  = X^3 +  6X^2 +  16X$ \\      %
36 & 0  &  0 & $Y^2         +   6Y  = X^3 +  6X^2 +  12X$ \\      %
49 & 1  &  0 & $Y^2 - 3XY           = X^3 +  3X^2 +   7X$ \\       %
\hline
\end{tabular}
\caption{Elliptic curves $E_N$ for $N$ with genus 1}
\label{table:4}
\end{table}

For each level $N$ with genus 1, we found that $(X,\,Y)=(G_{N,2},\,G_{N,3})$ %
is a point on an elliptic curve $E_N$ specified in Table~\ref{table:4} %
and verified up to $O(q^{20000})$. %
Moreover  $G_{N,M}=P_0(X)+P_1(X)Y$ where $P_0$ and $P_1$
are polynomials with degrees not exceeding $M/2$ and $(M-3)/2$,
respectively. 

With $N=21$, relations~(\ref{e21a},\ref{e21b},\ref{e21c}) show %
that $X=G_{21,2}$ is determined by an eta quotient and that $Y=G_{21,3}$ %
is determined by 3 eta quotients. The transformation $(X,\,Y)=(x-5,\,y-x-3)$ %
yields a minimal model $y^2+xy=x^3-4x-1$, whose small coefficients were noted in~\cite{NE}. %

Applying the {\tt ellak} procedure of {\tt Pari-GP} to $E_N$,
we reproduce the Fourier coefficients of $f_N=\sum_{M>0}C_{N,M}q^M$.
Thanks to work recorded at OEIS, we are able to provide formulas for $f_N$ in the 4 cases where a single
eta quotient does not suffice, namely for $N=17,19,21,49$:
\begin{gather}
f_{17}=\eta_1\eta_{17}\left(\psi_2\phi_{17}-\psi_{34}\phi_1\right),\quad
f_{19}=\left(\psi_4\phi_{38}-\psi_1\psi_{19}+\psi_{76}\phi_2\right)^2,\label{lat1}\\
f_{21}=\frac{3(\et73\et93+\et1{}\et72\et92\et{63}{}+\et12\et7{}\et9{}\et{63}2)}{2\et3{}\et{21}{}}
-\frac{\et34\et72}{2\et12}
+\frac{7\et3{}\et7{}\et{21}{3}}{2\et1{}}
-\frac{3\et34\et7{}\et{63}{}}{2\et1{}\et9{}},\label{lat2}\\
f_{49}=\theta_{7,14}^3\left(q\theta_{21,28}+q^2\theta_{14,35}-q^4\theta_{7,42}\right),\label{lat3}\\
\psi_n=\frac{\et{2n}2}{\et{n}{}},\quad \phi_n=\frac{\et{2n}{5}}{\et{n}2\et{4n}2},\quad
\theta_{a,b}=\sum_{n=-\infty}^\infty\left(-q^a\right)^{(n^2+n)/2}\left(-q^b\right)^{(n^2-n)/2},
\end{gather}
with~(\ref{lat3}) recorded in~\cite{BS}.
At $N=49$, we have complex multiplication, with $C_{49,p}=0$ for prime $p=3,5,6$ mod 7.
Moreover we have a pair of eta quotients,
\begin{equation}
x=\frac{\eta_{49}}{\eta_1{}}=\frac{G_{49,2}}{7},\quad
y=\frac{\et74}{\et14}=\frac{G_{49,7}}{7^2},
\end{equation}
with Fourier coefficients given by Rademacher sums. 
The latter is determined by $G_{49,2}$ and $G_{49,3}$. Hence the elliptic curve $E_{49}$ 
provides an algebraic relation between these eta quotients, namely
\begin{equation}
(2y-7x-35x^2-49x^3)^2=(4x+21x^2+28x^3)(1+7x+7x^2)^2.
\end{equation}

At $N=21$, we found 2937 eta quotients whose Fourier coefficients
are linear combinations of $R_{21,M}(n)$ with $M\le50$.  Including
the unit quotient, the tally of 2938 is the coefficient
of $x^{50}$ in the generating function
\begin{equation}
T_{21}(x)=\frac{1-x+x^2-x^3+2x^4}{(1-x)^2(1-x^4)^2}
\end{equation}
which predicts a total of 22126 eta quotients 
with Fourier coefficients determined  by $R_{21,M}(n)$ 
for $M\le100$. We have identified all of these.

At $N=36$, with 9 divisors, the corresponding tallies of eta quotients
are spectacularly large. Using Pad\'e approximants, the generating 
function was found to be 
\begin{gather}T_{36}(x)=\frac{H(x)+x^{18}H(1/x)}
{(1-x)^4(1-x^3)^2(1-x^4)^2(1-x^{12})},\\
H(x)=1-3x+6x^2-3x^3+6x^4+x^5+x^6+4x^7+4x^8+x^9,
\end{gather}
giving 49307076 eta quotients 
with Fourier coefficients determined by $R_{36,M}(n)$ 
for $M\le50$ and 8204657877 for $M\le100$. We were able to
identify all of these, by taking products of 78 eta quotients
found at $M\le12$ and eliminating redundancies.
Using more refined methods, we
also validated the monstrous tally of 180919436828
for $M\le150$.

From the denominator of $T_{36}(x)=\sum_{n\ge0}(c(M)+1)x^M$ it is clear
that the number $c(M)$ of non-trivial eta quotients determined by our procedures
may be found by polynomial interpolation at integers with the
same residue modulo 12. We denote the result by
$c(M)=p(M)+q_r(M)/12$ for $M = r$ mod 12, with a leading polynomial
\begin{equation}
p(n)=\frac{(n+1)(n+3)(n+5)(n+7)(n+9)(n+11)((n+6)^2-7)}{1935360}-1
\end{equation}
of degree 8 and sub-dominant terms that are at most quadratic:  
\begin{eqnarray}
q_0(n)&=&q_3(n)+q_4(n)+5,\\
q_1(n)&=&q_5(n)=q_7(n)=q_{11}(n)=0,\\
q_2(n)&=&q_{10}(n)=\frac{2(n+6)^2-5}{512},\\
q_3(n)&=&q_9(n)=\frac{(n+3)(n+9)}{9},\\
q_4(n)&=&q_8(n)=q_2(n)+\frac{(n+4)(n+8)}{16},\\
q_6(n)&=&q_2(n)+q_3(n)+1.
\end{eqnarray}

The situation at the prime levels $N=11,17,19$ is rather different.
Here we have a wealth of Rademacher sums but only one eta quotient.
Consider the case $N=19$. Since $f_{19}=q+O(q^3)$, we have $C_{19,2}=0$
and hence $R_{19,2}(n)$ yields integers. As noted, the sequence
$R_{19,2}(n)/19$ begins with
$$1, 4, 10, 25, 55, 116, 229, 440, 809, 1455, 2541, 4354, 7300, 12050,\ldots$$
This sequence may be developed using $\et{19}4/\et14$, 
which is determined by $G_{19,2}$ and $G_{19,3}$.
The elliptic curve relating the latter pair gives an algebraic relation between
$G_{19,2}$ and the eta quotient, namely
\begin{equation}
s^3/e_{19}=1+8s+19e_{19},\quad s=G_{19,2}/19,\quad e_{19}=\et{19}4/\et14,
\end{equation}
from which the expansion of $s=q+O(q^2)$ is easy developed, iteratively.

Similarly, at $N=11$ and $N=17$ we obtain the algebraic relations
\begin{eqnarray}
t^5/e_{11} &=& 1+13t+34t^2+11^2e_{11},\quad  t=G_{11,2}/11^2,\quad e_{11}=\et{11}{12}/\et1{12},\\
u^4/e_{17}&=&16+64u+34u^2-17^2e_{17},\quad u=G_{17,2}/17,\quad   e_{17}=\et{17}6/\et16,
\end{eqnarray}
and hence develop the expansions of $t=q+O(q^2)$ and $u=2q+O(q^2)$.

Intermediate between the plethora of eta quotients at $N=36$, with 9 divisors, and
their relative scarcity at $N=11,17,19,49$, with less than 4 divisors, sit the remaining
genus 1 levels, $N=14,15,20,21,24,27,32$. Proceeding as for $N=21$, we found the
generating functions 
\begin{eqnarray}
T_{14}(x)&=&\frac{1+x^3+2x^4+x^6+x^7}{(1-x)(1-x^2)(1-x^3)(1-x^6)},\\ 
T_{15}(x)&=&T_{21}(x)=\frac{1-x+x^2-x^3+2x^4}{(1-x)^2(1-x^4)^2},\\
T_{20}(x)&=&\frac{1-x+x^2+4x^3+2x^4+3x^6+x^7+x^8}{(1-x)^2(1-x^2)^2(1-x^3)(1-x^6)},\\
T_{24}(x)&=&\frac{(1+x^3)\left(h(x)+x^6h(1/x)\right)}{(1-x)^3(1-x^2)^3(1-x^4)^2},\quad
h(x)=1-2x+3x^2+2x^3,\\
T_{27}(x)&=&\frac{1-x^{11}}{(1-x)(1-x^2)(1-x^3)(1-x^5)(1-x^6)},\\
T_{32}(x)&=&\frac{1-x+x^2+2x^3+x^4-x^5+x^6}{(1-x)^2(1-x^2)^2(1-x^4)^2},
\end{eqnarray}
with the coefficient of $x^m$ in $T_{N}(x)$ giving the number of
eta quotients whose Fourier coefficients are determined by linear
combinations of the Rademacher sums $R_{N,M}(n)$ with $M\le m$.

\subsection{Rational Rademacher sums}

There are 5 levels with genus greater than 0 for which it appears that
the Rademacher sums $R_{N,M}(n)$ are rational for all positive integers
$M$ and $n$, namely $N=27,32,36,49$ with genus 1 and $N=64$ with genus 3. 
At genus 1, we convert these rationals to the integers
\begin{equation}
\overline{R}_{N,M}(n)=R_{N,M}(n)-C_{N,M}R_{N,1}(n),
\end{equation}
which vanish at $M=1$. 
The rationals $R_{64,M}(n)$ do not form integer sequences for $M=1,2,5$~mod~8.
To remedy this, we define
\begin{gather}
g_{8k+r}=\sum_{n>0}\left(R_{64,8k+r}(n)-c_{k,r}R_{64,r}(n)\right)q^n,\\
c_{k,1}=C_{32,8k+1},\quad c_{k,2}=C_{32,4k+1},\quad c_{k,5}=-C_{32,8k+5}/2,
\end{gather}
with $k\ge0$, $r\in[1,8]$ and $c_{k,r}=0$ for $r=3,4,6,7,8$. 
Then $g_M$ has integer Fourier coefficients, which vanish for $M=1,2,5$. Eta quotients appear in
\begin{equation}\frac{g_3}{2}=\frac{\et2{}\et42}{\et12\et8{}}-1,\quad
\frac{g_4}{8}=\frac{\et2{}\et{16}2}{\et12\et8{}},\quad
\frac{g_6}{4}=\frac{\et22\et44}{\et14\et82}-1,\quad
\frac{g_7}{4}=\frac{\et27}{\et16\et8{}}-1,\quad
\frac{g_8}{32}=\frac{\et22\et84}{\et14\et42}.
\end{equation}
Moreover, $g_6=(4+g_3)g_3,~g_7=2g_3+(2+g_3)g_4,~g_8=(4+g_4)g_4$.

We found that $(X,Y)=(g_3+2,g_4+2)$ is a point on the genus 3 curve
\begin{equation}
Y(Y^2+4)=X^4\label{N64}
\end{equation}
and that $g_M=P_0(X)+P_1(X)Y+P_2(X)Y^2$, with $P_n$ a polynomial
of degree at most $\lfloor(M-4n)/3\rfloor$.
Every eta quotient of that form has Fourier coefficients that are linear combinations
of Rademacher sums at $N=64$. 

\subsection{Genus 2}

It is instructive to compare the genus 3 case $N=64$
with the genus 2 case $N=50$. For the latter  
we construct integer sequences as follows:
\begin{gather}\overline{R}_{50,M}(n)=\left\{\begin{array}{l}
R_{50,M}(n)-d(M)R_{50,1}(n),~{\rm for}~M=1,4~{\rm mod}~5,\\
R_{50,M}(n)+d(M)R_{50,2}(n),~{\rm for}~M=2,3~{\rm mod}~5,\\
R_{50,M}(n),~{\rm for}~M=0~{\rm mod}~5,\end{array}\right.\\
\sum_{M>0}d(M)q^M=f_{50}=
q - q^2 + q^3 + q^4 - q^6 + 2q^7 - q^8 - 2q^9 + O(q^{11}),
\end{gather}
where $f_{50}$ is the weight 2 level 50 Hecke eigenform whose Fourier
coefficients $d(M)$ are obtained from the L-series of the elliptic curve
$y(y+x+1)=x^3-x-2$. Then $G_{50,M}=\sum_{n>0}\overline{R}_{50,M}(n)q^n$
vanishes by construction at $M=1,2$ and yields eta quotients at $M=3,5$, with
\begin{equation}
\frac{G_{50,3}}{5}=\frac{\et2{}\et{25}2}{\et12\et{50}{}}-1,\quad
\frac{G_{50,5}}{20}=\frac{\et2{}\et{10}3}{\et13\et5{}}.
\end{equation}
The Fourier coefficients of $G_{50,4}$
are also identified by an eta quotient: $\overline{R}_{50,4}(n)/10$ is the
coefficient of $q^{2n}$ in the Fourier expansion of $\eta_{25}/\eta_1$.

We found that $(X,Y)=(G_{50,3},G_{50,4})$ is a point on the curve
\begin{equation}
Y^3+4(X+5)Y^2+2(X+5)(X+10)Y=X(X+5)(X^2+8X+20)
\end{equation}
which {\tt Sage} confirmed as having genus 2.

Proceeding similarly for $(X,Y)=(G_{N,3},G_{N,4})$ we found the curves
\begin{gather}
Y^3+55Y^2-2(X^2+11X-484)Y=X(X^3+34X^2+473X+2904)\\
Y^3+2(4X+69)Y^2+(9X^2+460X+4761)Y=X(X^3+55X^2+1035X+6348)\\
Y^3+8(X+13)Y^2+4(X+13)(3X+52)Y=X(X+13)(X^2+28X+208)\\
Y^3+21Y^2+(5X^2+70X+392)Y=X(X+14)(X^2+14X+56)\\
Y^3+2(2X+29)Y^2-(4X^2-58X-841)Y=X(X^3+31X^2+406X+1682)\\
Y^3+2(4X+31)Y^2+(11X^2+217X+961)Y=X(X^3+31X^2+310X+961)\\
Y^3-8XY^2-X(2X+259)Y=X(X^3-10X^2+148X+1369)
\end{gather}
at $N=22,23,26,28,29,31,37$, respectively. All have genus 2.

\subsection{Genus 3}

We found these genus 3 curves for
$N=30,33,34,35,39,40,41,43,45,48,64$:
\begin{gather}
Y^2(Y-2X)(Y-3X)+X(7X^2-30X+75)Y\nonumber\\
=X^5+25X(4X-5),~{\rm with}~(X,Y)=(G_{30,4}+5,G_{30,5}+10),\\
Y^4+X(5X-11)Y^2-X^2(4X-11)Y\nonumber\\
=X^3(X^2-11X+22),~{\rm with}~(X,Y)=(G_{33,4}+11,G_{33,5}+11),\\
Y^4+10XY^3+X(21X-221)Y^2+2X(3X^2-119X+867)Y\nonumber\\
=X(X^4-2X^3+51X^2-578X+4913),~{\rm with}~(X,Y)=(G_{34,4}+17,G_{34,5}+17),\\
Y^4+10(X-3)Y^3+(31X^2-210X+800)Y^2+(12X^3+25X^2+200X-2000)Y\nonumber\\
=X(X^4-5X^3-15X^2-200X-2000),~{\rm with}~(X,Y)=(G_{35,4}+20,G_{35,5}+10),\\
Y^4+5XY^3+3X(X+13)Y^2-X(19X^2-234X+507)Y\nonumber\\
=X(X^4-14X^3+234X^2-1690X+2197),~{\rm with}~(X,Y)=(G_{39,4}+13,G_{39,5}+13),\\ 
Y^4=X(X+5)(X^2(X+4)-4Y^2),~{\rm with}~(X,Y)=(G_{40,4},G_{40,5}),\\
Y^4+(10X-41)Y^3+X(30X-451)Y^2+X^2(11X-1681)Y\nonumber\\
=X^3(X^2+70X+2214),~{\rm with}~(X,Y)=(G_{41,4},G_{41,5}+41),\\
32Y^4+(40X+43)Y^3+(94X^2+1591X+9245)Y^2+X(49X^2+946X+5547)Y\nonumber\\ 
=X^3(X^2+21X+129),~{\rm with}~(X,Y)=(G_{43,3},G_{43,5}-2G_{43,3}),\\        
(Y^2+5X)^2=X^3(X^2-Y),~{\rm with}~(X,Y)=(G_{45,4}+5,G_{45,5}+5),\\
Y^4=X^3(X-3)(X-4),~{\rm with}~(X,Y)=(G_{48,4}+6,G_{48,5}+6),\\
Y(Y^2+4)=X^4,~{\rm with}~(X,Y)=(G_{64,3}+2,G_{64,4}+2),
\end{gather}
where the final curve at $N=64$ was already given in~(\ref{N64})
and was obtained by subtractions that make $G_{64,M}$ vanish at $M=1,2,5$.
At $N=43$, the subtractions make $G_{43,M}$ vanish at $M=1,2,4$. In all other
cases with genus 3, $G_{N,M}$ vanishes for $M=1,2,3$.

\subsection{Genus 4}

At $N=81$, we found that $(X,Y)=(G_{81,5},G_{81,6})$ lies on the genus 4 curve
\begin{equation}
Y^3(Y+3)^2+3(X+3)^3(Y+3)(2Y+3)=(X+3)^6,\quad \frac{Y}{3}=\frac{\et34}{\et13\et9{}}-1.
\end{equation}
The Fourier coefficients of $X/9$ form an integer sequence beginning with
$$1, 2, 4, 7, 13, 21, 35, 55, 87, 132, 200, 295, 434, 625, 897, 1267, 1782, 2475,\ldots$$
for $n=1$ to 18. The general term is given by Rademacher sums as
\begin{equation}
\frac{R_{81,5}(n)+R_{81,2}(n)}{9}=\frac{\exp(4\pi\sqrt{5n}/9)}{27(4n^3/5)^{1/4}}
\left(1-\frac{27}{32\pi\sqrt{5n}}+O(1/n)\right).
\end{equation}
The Fourier coefficients of $Y/9$ form an integer sequence beginning with
$$1, 3, 6, 13, 24, 45, 77, 132, 216, 351, 552, 861, 1313, 1986, 2952, 4354, 6336,\ldots$$
for $n=1$ to 17. The general term is given by a Rademacher sum
\begin{equation}
\frac{R_{81,6}(n)}{9}
=\frac{\exp(4\pi\sqrt{6n}/9)}{27(2n^3/3)^{1/4}}
\left(1-\frac{27}{32\pi\sqrt{6n}}+O(1/n)\right).
\end{equation}

\subsection{Genus 5}

At $N=72$, with genus 5, we obtain integer Fourier coefficients in
\begin{gather}
G_{72,M}=\sum_{n>0}\overline{R}_{72,M}(n)q^n,\quad
\overline{R}_{72,M}(n)=R_{72,M}(n)-\sum_{r=1,2,3,5,7}p_r(M)R_{72,r}(n),\\
\et{12}2\et{18}4/\et{36}2=\sum_{n>0}p_1(n)q^n,\quad
\et{12}{4}=\sum_{n>0}p_2(n)q^n,\quad
\eta_6\eta_{12}\eta_{18}\eta_{36}=\sum_{n>0}p_3(n)q^n,\\
\et62\et{36}4/\et{18}2=\sum_{n>0}p_5(n)q^n,\quad
\et{6}{4}=\sum_{n>0}\left(p_1(n)-4p_7(n)\right)q^n.
\end{gather}
Then $G_{72,M}$ vanishes for $M=1,2,3,5,7$. 
Moreover $(G_{72,4},\,G_{72,6})$ is a point on the elliptic 
curve $E_{36}$, while $(G_{72,6},\,G_{72,9})$ lies on $E_{24}$. 
Eliminating $G_{72,6}$, we obtain a genus 5 curve from the resultant:
\begin{gather}
\left(Y^2+12Y+34-2(X+2)^3\right)^2=\left((X+2)^3+1\right)\left((X+2)^3-2\right)^2,\label{G5}\\
X=G_{72,4}=\frac{6\et2{}\et3{}\et{18}2}{\et12\et6{}\et9{}},\quad
Y=G_{72,9}=\frac{6\et25\et32\et6{}}{\et16\et4{}\et{12}{}}-6.
\end{gather}
The Fourier coefficients of $X/6$ are the integers $R_{18,1}(n)/6={\rm A}128129(n)$.
The Fourier coefficients of $Y/12$ form an integer sequence beginning with
$$3, 11, 33, 87, 210, 473, 1008, 2055, 4035, 7674, 14196, 25629, 45282, 78472,\ldots$$
for $n=1$ to 14. The general term is given by Rademacher sums as
\begin{equation}
\frac{R_{24,3}(n)+R_{24,1}(n)}{12}=\frac{\exp(\pi\sqrt{2n})}{24(2n^3)^{1/4}}
\left(1-\frac{3}{8\pi\sqrt{2n}}+O(1/n)\right).
\end{equation}

\subsection{Genus 6}

Moving on to the genus 6 case $N=121$, we determined that subtractions
of $R_{121,r}(n)$ are needed for the 6 values $r=1,2,3,4,6,11$.
The coefficients of these subtractions are determined by 
four new forms and two old forms of weight 2 and level 121.
The new forms are the L-series of the elliptic curves
$y^2+xy+y=x^3+x^2-30x-76$, $y^2+y=x^3-x^2-7x+10$,
$y^2+xy=x^3+x^2-2x-7$ and $y^2+y=x^3-x^2-40x-221$.
The old forms are $(\eta_1\eta_{11})^2$ and
$(\eta_{11}\eta_{121})^2$. The first two non-zero integer series are
\begin{equation}
\overline{R}_{121,5}(n)=R_{121,5}(n)-R_{121,4}(n)-R_{121,3}(n),\quad
\overline{R}_{121,7}(n)=R_{121,7}(n)+R_{121,6}(n).
\end{equation}
We expect their generators, $G_{121,5}$ and $G_{121,7}$, to define a curve of genus 6
with degree 7 in $G_{121,5}$ and degree 5 in $G_{121,7}$. This is indeed the case.
We found that $(X,Y)=(G_{121,5}/11,G_{121,7}/11)$ is a point on the curve
\begin{gather}
Y^5-20XY^4+5X(10X-9)Y^3=X(132X^3-64X^2-33X+31)Y^2\nonumber\\
+X(33X^4+95X^3-48X^2-5X+9)Y\nonumber\\
+X(121X^6-66X^5+23X^4+18X^3-9X^2+1)
\end{gather}
which {\tt Sage} confirmed as having genus 6. Since integer
combinations of Rademacher sums are computable with ease, 
we are able to validate this curve up to $O(q^{1000})$ in a matter of seconds. 

\subsection{Genus 7}

In the genus 7 case $N=100$, we determined that subtractions
of $R_{100,r}(n)$ are required for the 7 values $r=1,2,3,4,5,7,9$,
with  coefficients determined by 
6 old forms of weight 2 and level 100 and a new form
\begin{equation}
f_{100}=q+2q^3-2q^7+q^9-2q^{13}+6q^{17}-4q^{19}-4q^{21}-6q^{23}-4q^{27}+6q^{29}+O(q^{31})
\end{equation}
which is the L-series of the elliptic curve $y^2=x^3-x^2-33x+62$.

We found that $(X,Y)=(G_{100,6}+5,G_{100,15}+10)$ lies on the genus 7 curve 
\begin{gather}
Y^6=X(3X^4-90X^3+415X^2-560X+200)Y^4\nonumber\\
-X(3X^9+20X^8-350X^7+1795X^6-4790X^5\nonumber\\
+7805X^4-8350X^3+7325X^2-4625X+1375)Y^2\nonumber\\
+X(X^2+2X+5)(X^4-6X^3+14X^2-10X+5)^2(X^4-5X^3+15X^2-25X+25),
\end{gather}
with $X/5=\et2{}\et{25}2/(\et12\et{50}{})$ and the Fourier coefficients of $Y/10$
given by
$$1, 6, 26, 88, 258, 686, 1688, 3904, 8594, 18142, 36946, 72952, 140184, 262948,\ldots$$
for $n=0$ to 13, and in general by $(R_{100,15}(n)+2R_{100,5}(n))/10$ for $n>0$. 

\subsection{Genus 8}

When $N=p^2$ with prime $p=12k+1$, the genus of $\Gamma_0(N)$ is given 
by $g_0(N)=3k(4k-1)-1$. Setting $k=1$, we obtain $g_0(169)=8$. Moreover $N=169$
is the largest level with genus 8.
We devised a procedure of 8 subtractions that gives
integer sequences $\overline{R}_{169,M}(n)$
by subtracting multiples of $R_{169,r}(n)$,
with $r=1,2,3,4,5,6,8,9$.
The subtraction coefficients are determined
by 8 modular forms of level 169 and weight 2,
of which two have Fourier coefficients in ${\bf Q}({\sqrt3})$.
The rest have coefficients in the
cubic number fields $x(x^2-1)=\pm(1-2x^2)$.

Our first non-zero integer sequence occurs at $M=7$, where
\begin{equation}
\overline{R}_{169,7}(n)=R_{169, 7}(n)-R_{169,6}(n)-R_{169,5}(n)+R_{169,2}(n)
\end{equation}
is the coefficient of $q^n$ in $G_{169,7}=13\eta_{169}/\eta_1.$ 
There is no subtraction at $M=13$, where 
$G_{169,13}=G_{13,1}=13\et{13}2/\et12.$  

We found that $(X,Y)=(G_{169,7}/13,G_{169,13}/13)$ is point on a
genus 8 curve with degree 13 in $X$ and degree 7 in $Y$, namely
\begin{gather}
Y^7=143XY^6+156X(39X^2-17X+3)Y^5+X\sum_{k=2}^6P_k(X)Y^{6-k},\\
P_2(X)=26(3211X^4-2249X^3+819X^2-173X+19),\\
P_3(X)=26(21970X^6-18759X^5+8619X^4-2743X^3+663X^2-111X+10),\\
P_4(X)=26(169X^4-104X^3+39X^2-8X+1)(507X^4-169X^3+26X^2-13X+3),\\
P_5(X)=13(371293X^{10}-371293X^9+199927X^8-81289X^7\nonumber\\
+28561X^6-8619X^5+2197X^4-481X^3+91X^2-13X+1),\\
P_6(X)=4826809X^{12}-4826809X^{11}+2599051X^{10}-1113879X^9+428415X^8\nonumber\\
-142805X^7+41743X^6-10985X^5+2535X^4-507X^3+91X^2-13X+1.
\end{gather}

\subsection{Genus 13}

Finally, we studied $N=144$, with genus 13. Integer sequences are obtained
after subtractions at $r=1, 2, 3, 4, 5, 6, 7, 9, 10, 11, 13, 14, 19$, 
with coefficients determined by 11 old forms and two new forms.
The old forms were already encountered as eta quotients at $N=48$ and $N=72$.
They determine the subtractions for $R_{144,M}(n)$ when $M$ is divisible by 3 or 2.
The new forms relate to the subtractions that are needed when $M$ is coprime to 6.
One of these new forms is the eta quotient 
\begin{equation}
f_{144a}=\frac{\et{12}{12}}{\et46\et{24}4}=q + 4q^7 + 2q^{13} - 8q^{19} - 5q^{25} + 4q^{31} + \ldots
\label{f144a}
\end{equation}
The other has Fourier coefficients determined by the L-series of the elliptic curve
$y^2=x^3+6x+7$, which gives 
\begin{equation}
f_{144b}=q + 2q^5 + 4q^{11} - 2q^{13} - 2q^{17} + 4q^{19} - 8q^{23} - q^{25} - 6q^{29} - 8q^{31} + \ldots
\end{equation}
We reduced this to a telling combination of 5 eta quotients with weight 2:
\begin{equation}
f_{144b}=
\frac{\et{24}3\et{36}2\et{72}{}}{\et{48}{}\et{144}{}}
+\frac{2\et{12}{}\et{24}2\et{72}4}{\et{36}{}\et{48}{}\et{144}{}}
+\frac{4\et{12}{2}\et{48}{}\et{72}{}\et{144}{}}{\et{24}{}}
-\frac{2\et{12}{}\et{24}3\et{72}{}\et{144}2}{\et{36}{}\et{48}2}
+\frac{4\et{12}2\et{144}4}{\et{72}2}\label{f144b}
\end{equation}
with $C_1(M)+2C_5(M)+4C_{11}(M)-2C_{13}(M)+4C_{19}(M)$
giving the coefficient of $q^M$ in~(\ref{f144b}),  
where $C_r(M)R_{144,r}(n)$ is subtracted from $R_{144,M}(n)$
to make  $\overline{R}_{144,M}(n)$ an integer sequence.
The subtraction at $r=7$ is determined by the eta quotient in~(\ref{f144a})
where the coefficient of $q^M$
is $C_1(M)+4C_7(M)+2C_{13}(M)-8C_{19}(M)$.
Hence we determine all the subtractions by eta quotients.

We are left with 6 values of $M<20$ for which
$G_{144,M}=\sum_{n>0}\overline{R}_{144,M}(n)q^n$
is non-zero, namely $M=8,12,15,16,17,18$.  
To produce a genus 13 curve we should choose a coprime
pair of $M$ values. The simplest choice is the pair $(8,15)$. 
We know that $(G_{144,8},G_{144,18})=(G_{72,4},G_{72,9})$
gives a point on the genus 5 curve~(\ref{G5}) found at $N=72$.
Moreover $(G_{144,15},G_{144,18})=(G_{48,5},G_{48,6})$ gives a point
on a genus 3 curve that is not hard to determine. Then, by taking a resultant
to eliminate $G_{144,18}$, we determined that
$(X,Y)=(G_{144,8}+2,G_{144,15}+6)$
lies on the genus 13 curve
\begin{equation}
(Y^4-8(X^3+1)^2)^2=(X^3+1)(X^6+20X^3-8)^2
\end{equation}
neatly parametrized by eta quotients as follows
\begin{equation}\frac{X}{2}=\frac{\et23\et3{}}{\et13\et6{}},\quad 
\frac{Y}{6}=\frac{\et22\et64}{\et14\et{12}2}.
\end{equation}

\subsection{Remarks} %

{\bf Remark 1.} 
After we completed this work, Yajun Zhou kindly called our             %
attention to~\cite[Theorem 8.12]{Moonshine}. The methods            
used there may be capable of furnishing proofs of some of              %
our empirical findings in Section 3, following the approach              %
that Knopp~\cite{K} attributes to Rademacher~\cite{R1} as            %
an ``entirely fresh viewpoint",                                                      %
namely by adopting formulas~(\ref{RNM},\ref{Kc})  as                   %
definitions of Fourier coeffcients of objects $G_{N,M}$ and               %
demonstrating that the latter have the required modular properties.  %
At genus 0, with a unique Hauptmodul, that could furnish a proof     %
of Table 1. At higher genera, more work might be needed.                %

\noindent{\bf Remark 2.}
We conclude this section with a note on the approach in~\cite{YY1,YY2}
to modular curves. In~\cite[Section~4.1]{YY2}, Yifan Yang gives
modular curves, up to level $N=50$, that are parametrized by quotients of ``generalized'' Dedekind
eta functions~\cite{YY1}, in the many cases where the eta function itself is insufficient to solve the problem.
Moreover his $q$-expansions are highly singular as $q\to0$. Our approach was quite different. We began
with an explicit formula~(\ref{RNM}) that reproduces, at $M=1$, the Fourier coefficients of the genus 0 eta quotients
taken as ``canonical'' Hauptmoduln in~\cite[Table~8]{Maier}, which vanish as $q\to0$. At genus 1,
after subtraction of the non-integer sequence $R_{N,1}(n)$, we obtained $G_{N,M}=\sum_{n>0}\overline{R}_{N,M}(n)q^n$ 
as Fourier series with integer coefficients, vanishing at $q=0$. Then $G_{N,2}$ and $G_{N,3}$ parametrize our modular curve.
We were able to extend this to higher genera. It may be that our explicit Fourier coefficients
are capable of reproducing those of Yang's ``generalized'' Dedekind eta quotients, after performing a
Fricke involution $z\longmapsto-1/(Nz)$ on his Ans\"atze. We have not investigated this, since it lay outside the remit of our title.
  
\section{Conclusions}

\begin{enumerate}
\item Eta quotients on $\Gamma_0(6)$, with 4 cusps, neatly solve the equal-mass two and three loop sunrise problems,
whose differential equations with respect to the external energy have 4 singular points. This cannot continue, since
at higher loops there is more than one pseudo-threshold.
\item Atkin-Lehner transformations on $\Gamma_0(6)$ yield optimal nomes.
\item For the on-shell problem, Eichler integrals of eta quotients on $\Gamma_0(6)$ yield Bessel moments at 4 and 6 loops
that are periods or quasi-periods. 
\item Rademacher sums yield the Fourier coefficients of a Hauptmodul for $\Gamma_0(6)$ and for all other levels of genus 0.
\item After subtractions determined by weight 2 cusp forms, they yield the Fourier coefficients of vast numbers of eta quotients.
\item They yield the Fourier coefficients of parametrizations of modular curves, irrespective of whether the Fourier series
are eta quotients. 
\end{enumerate}

\noindent{\bf Acknowledgements.}
The second author thanks KMPB for hospitality and colleagues at conferences
in Zeuthen, Bonn, St. Goar and Les Houches for advice and encouragement that emboldened
our joint effort to tackle eta quotients beyond the remit of genus zero so far encountered in 
massive Feynman diagrams. We especially thank Johannes Bl\"umlein for his question 
on the possibility of obtaining an explicit formula for Fourier coefficients of the Hauptmodul of
$\Gamma_0(6)$ and Freeman Dyson for urging us to try to emulate the notable work by 
Rademacher on partition numbers~\cite{R2}. 
We thank Yajun Zhou and an anonymous referee for helpful suggestions that improved our presentation. %

\end{document}